%% file: GordinaHaga2.tex
\documentclass
{amsart}%
\usepackage[hypertex]{hyperref}
\usepackage{cite}
\usepackage{amssymb}
\usepackage{amsmath}
\usepackage{amsthm}
\usepackage{amsfonts}
\usepackage{mathrsfs}
\usepackage{graphicx}%
\usepackage{color}
\usepackage{eucal,amscd}
\usepackage[toc,page]{appendix}

\newtheorem{thm}{Theorem}[section]

\newtheorem{prop}[thm]{Proposition}

\theoremstyle{definition}

\newtheorem{df}[thm]{Definition}

\newtheorem{example}{Example}[section]

\numberwithin{equation}{section}

\begin{document}

\input{GordinaHaga-title}

\section{Introduction}

Let $K$ be a Lie group with the identity $e$ and let $X_{t}$ be a L\'evy process with values in $K$ starting at $e$.  One obtains a semigroup of operators $(T(t),t\geqslant 0)$ on the Banach space $C_0(K)$ of functions on $K$ which vanish at infinity, by defining
\[
(T(t)f)(k)=\mathbf E(f(k X_{t}))
\]
for each $t\geqslant 0$, $k \in K$ and $f\in C_0(K)$.

When $K=\mathbb{R}^n$ the characteristic function of the process $X_{t}$ is given by the L\'evy-Khintchine formula
\begin{align*}
\mathbf E(e^{iu\cdot X_{t}})=e^{t\varphi(u)}
\end{align*}
for all $u\in \mathbb{R}^n$, $t\geqslant 0$, where
\begin{align}
\varphi(u)=im\cdot u-\frac{1}{2}u\cdot au+\int_{\mathbb{R}^n-\{0\}}\left(e^{iu\cdot y}-1-i\frac{u\cdot y}{1+|y|^2}\right)\nu(dy).\label{levykhintchineintro}
\end{align}
Here $m\in \mathbb{R}^n$, $a$ is a non-negative symmetric $n\times n$ matrix and $\nu$ is a L\'evy measure on $\mathbb{R}^n\backslash\{0\}$ (see \cite{dudleybook} for details).

Let the differential operator $\mathcal D$ be defined on $C_0(\mathbb{R}^n)$ by $\mathcal D=(\mathcal D_1,\dots,\mathcal D_n)$ with $\mathcal D_j=\frac{1}{i}\frac{\partial}{\partial x_j}$.  The generator $\mathcal A$ of this semigroup satisfies the relation
\begin{align}
\mathcal A= \varphi(\mathcal D)\label{levygeneratorinrn}
\end{align}
where $\varphi$ is as in (\ref{levykhintchineintro}).  Indeed, one makes the observation that $\mathcal A$ is in fact a {\it pseudo-differential} operator (see \cite[pp. 139-170]{heilharmonicanalysis}) with symbol $\varphi(u)$ \cite{applebaumbook}.

There has been interest in extending this characterization of L\'evy processes to  Lie groups, including the book by M. Liao \cite{LiaoBook}.  For an arbitrary Lie group $K$, one may define the group Fourier transform $\hat f$ of a suitably chosen function by

\begin{align}
\hat f(\pi)=\int_K f(k)\pi(k)\operatorname{d} k, \label{grouptransform}
\end{align}
where $\pi$ is a unitary irreducible representation of $K$ and $dk$ is Haar measure (see \cite{FollandHABook}).  This Fourier transform may be inverted if a complete set of unitary irreducible representations of $K$ is known.  Fourier inversion is necessary when formulating a theory of pseudo-differential operators.  Because the representation theory of Lie groups is only fully understood for specific subclasses of Lie groups, \eqref{grouptransform} might have to be adopted to each case separately.

In the current paper we consider a step 3 nilpotent group. If $K$ is a general nilpotent group, then Kirillov's method of co-adjoint orbits provides explicit formulae of all unitary irreducible representations of $K$ (refer to Theorem \ref{kirillov} below for details).  In \cite{BeltitaBeltita2009b}, Belti\c ta and Belti\c ta apply this technique to describe the Weyl functional calculus for arbitrary nilpotent Lie groups.  In what follows, we make use of this symbolic calculus, and of general results from \cite{ApplebaumCohen2004} to describe the quantized generator $\mathcal L^\pi$ of a L\'evy processes $X_{t}$ in a step 3 nilpotent Lie group $G$.  We describe $G$ and the collection of all unitary irreducible representations of $G$ in Sections \ref{s2} and \ref{s3}.  In Section \ref{s4} we describe the Weyl functional calculus for $G$ and in Section \ref{s5} we prove the following theorem which is the main result of this paper.

\begin{thm}  The operator $\mathcal L^\pi$ is a pseudo-differential operator.  Moreover, the space $C^\infty_c(\mathbb{R})$ is a core for $\mathcal L^\pi$.
\end{thm}

Here $C^\infty _c(\mathbb{R})$ denotes the collection of infinitely differentiable functions of compact support on $\mathbb{R}$.

The case when $K$ is the Heisenberg group was treated in \cite{ApplebaumCohen2004}. Then one may make use of the classical Schr\"odinger representations.  The resulting pseudo-differential calculus (referred to as the classical Weyl functional calculus) has been used to express the generators of these semigroups as
\begin{align*}
(\mathcal A_\pi f)(x)=(2\pi )^{-n}\int_{\mathbb{R}^{2n}}\sigma\left(\frac{1}{2}(x+y),\xi\right)e^{i(x-y)\cdot \xi}f(y)dyd\xi.
\end{align*}
Here one works through the Schr\"odinger representations and deals not with the original semigroup generator, but with their images $\mathcal A_\pi$ (henceforth referred to as the quantization of the generator $\mathcal A$).  L\'evy processes in the Heisenberg group have been thoroughly investigated;  the reader may refer to the work of D. Applebaum and S. Cohen in \cite{ApplebaumCohen2004} for a complete treatment of the Heisenberg group case. We restricted our study to finding an explicit form of  the quantized generator, and as a result developed the method which is potentially applicable to a larger class of nilpotent groups.

\section{A Step 3 Nilpotent Lie Group}\label{s2}
Let $G$ denote $\mathbb{R}^4$ with the multiplication law
\begin{align*}
\{w_1,x_1,y_1,z_1\}*\{w_2,x_2,y_2,z_2\}=&\left\{w_1+w_2,\,x_1+x_2,\,y_1+y_2+w_1x_2,\, \phantom{\frac{x}{2}}\right.\\
&\left. z_1+z_2+w_1\left(y_2+\frac{w_1x_2}{2}\right) \right\}.
\end{align*}
With respect to this operation $G$ is a Lie group with identity $\{0,0,0,0\}$ and inversion given by

\[
\{w,x,y,z\}^{-1}=\left\{-w,-x,-y+wx,-z+w\left(y-\frac{wx}{2}\right)\right\}.
\]
The Lie algebra $\mathfrak{g}$ of left invariant vector fields of $G$ is spanned by $\{W,X,Y,Z\}$ where
\begin{align*}
&W=\frac{\partial}{\partial w},	
\\
&X=\frac{\partial}{\partial x}+w\frac{\partial}{\partial y}+\frac{w^2}{2}\frac{\partial}{\partial z},	
\\
&Y=\frac{\partial}{\partial y}+w\frac{\partial}{\partial z},
\\
&Z=\frac{\partial}{\partial z}.	
\end{align*}
These vector fields satisfy the following commutation relation
\begin{align*}
& [W,X]=Y,
\\
& [W,Y]=Z,
\end{align*}
with all other brackets zero.  This Lie algebra $\mathfrak{g}$ is step 3 nilpotent, and the exponential map $\exp:\mathfrak{g}\rightarrow G$ is given by
\[
\exp( w, x, y, z )=\left\{w, x, y+\frac{wx}{2}, z+\frac{xy}{2}+\frac{w^2x}{6}\right\}.
 \]
 Because both the underlying manifold of $G$ and $\mathfrak{g}$ are $\mathbb{R}^4$, we adopt the convention of $\{w,x,y,z\}$ when referring to a point in $G$ and $( w, x, y, z )$ when referring to a point in $\mathfrak{g}$.

Any Lie group naturally acts on its Lie algebra via the adjoint representation.  The adjoint action of $G$ on $\mathfrak{g}$ is given by

\begin{align*}
\operatorname{Ad}(\{w,x,y,z\})(a,b,c,d)=\left(a,\,b,\,c+(wb-ax),\, d+(wc-ay)+\frac{w^2b}{2}\right).
\end{align*}
Let $\mathfrak{g}^*$ denote the linear dual of $\mathfrak{g}$.  The adjoint action induces the {\it co-adjoint action of} $G$ {\it on} $\mathfrak{g}^*$, defined for each $l\in \mathfrak{g}^*$ as

\begin{align*}
\operatorname{Ad}^*(\{w,x,y,z\})\left(l(a,b,c,d)\right)=l(\operatorname{Ad}(\{w,x,y,z\}^{-1})(a,b,c,d)).
\end{align*}
In the following sections we will make use of unitary irreducible representations of $G$.  These representations for nilpotent Lie groups can be classified by using Kirillov's method of co-adjoint orbits, as stated in Theorem \ref{kirillov}.  To this end we need to identify the co-adjoint orbits of $G$ in $\mathfrak{g}^*$.  If
\[
l(a,b,c,d)=\alpha a+\beta b+\gamma c+\delta d
 \]
for $(a,b,c,d)\in \mathfrak{g}$ then we will adopt the convention of writing $l=[\alpha,\beta,\gamma,\delta]$.  In these coordinates, the co-adjoint action of $G$ on $\mathfrak{g}^*$ is given by
\begin{align}
& \operatorname{Ad}^*(\{w,x,y,z\})[\alpha, \beta, \gamma, \delta]=\label{caaction}
\\
&\left[\alpha+x\gamma+(y-wx)\delta, \beta-w\gamma+\frac{w^2\delta}{2}, \gamma-w\delta, \delta \right].\notag
\end{align}

\section{Representation Theory of $G$}\label{reps}\label{s3}
 To describe a complete set of unitary irreducible representations of $G$, we make use of the fact that $G$ is nilpotent.  The following result, due to Kirillov is presented in  \cite[Section 2.2]{CorwinGreenleafBook}.
\begin{thm}\label{kirillov}(Kirillov)
Let $K$ be any connected nilpotent Lie group with Lie algebra $\mathfrak k$.
\begin{enumerate}
\item If $l\in \mathfrak k^*$ then there exists a subalgebra $\mathfrak m_l$ of $\mathfrak k$ of maximal dimension such that $l([m_1,m_2])=0$ for all $m_1,m_2\in \mathfrak m_l$.
\item $M_l=\exp(\mathfrak m_l)$ is a closed subgroup of $K$, and $\rho_l(\exp(m))=e^{2\pi il(m)}$ is one dimensional representation of $M_l$.
\item The induced representation $\operatorname{Ind}_{M_l,\rho_l}^K$ is a unitary irreducible representation of $K$.
\item If $\pi$ is any unitary irreducible representation of $K$, then there exists $l\in \mathfrak k^*$ such that $\pi$ is unitarily equivalent to $\operatorname{Ind}_{M_l,\rho_l}^K$.
\item Two irreducible representations $\pi_1=\operatorname{Ind}_{M_{l_1}, \rho_{l_1}}^K$ and $\pi_2=\operatorname{Ind}_{M_{l_2}, \rho_{l_2}}^K$ are unitarily equivalent if and only if $l_1$ and $l_2$ are elements of the same coadjoint orbit of $K$ in $\mathfrak k^*$.
\end{enumerate}
\end{thm}
If $l$ and $\mathfrak m_l$ are as in Theorem \ref{kirillov}, then the subalgebra $\mathfrak m_l$ is said to be a {\it maximal subordinate algebra} for $l$.

Theorem \ref{kirillov} implies that the set of unitary irreducible representations of $G$ is indexed by the set of co-adjoint orbits of $G$ in $\mathfrak g^*$.  The coadjoint action described by (\ref{caaction}) allows for an explicit parametrization of these orbits.  This parametrization can be used to give an explicit expression of unitary dual of $G$, as presented in the following proposition.  This calculation can be found in \cite{CorwinGreenleafBook,KirillovBook}, but we include it here for completeness.

\begin{prop}\label{unirrepclassification}
If $\pi$ is a unitary irreducible representation of $G$, then $\pi$ is unitarily equivalent to a representation of one of the following classes.

\begin{enumerate}
  \item[]Class 1. $\pi$ is a unitary character of $G$ given by \[\pi(\{w,x,y,z\})(z)=e^{2\pi i(\alpha w + \beta x)}z\] for some $\alpha,\beta\in \mathbb{R}$ and any $z\in \mathbb{C}$.
  \item[]Class 2. $\pi$ is  a representation on $L^2(\mathbb{R})$ given by \[\pi(\{w,x,y,z\})f(k)=e^{2\pi i \gamma \left(y+\frac{kx}{2}\right)}f(k+w)\] for some $\gamma\in \mathbb{R}$.
  \item[]Class 3. $\pi$ is a representation on $L^2(\mathbb{R})$ given by \[\pi(\{w,x,y,z\})f(k)=e^{2\pi i  \left(\beta x + \delta \left(z+k\left(y+\frac{kx}{2}\right)\right)\right)}f(k+w)\] where $\delta\in \mathbb{R}^\times$, $\beta\in \mathbb{R}$.
\end{enumerate}

\end{prop}

\begin{proof}
If $[\alpha,\beta,\gamma,\delta]\in \mathfrak{g}^*$ and $\{w,x,y,z\}\in G$ then $\pi_{[\alpha,\beta,\gamma,\delta]}(\{w,x,y,z\})$ can be computed by considering some individual cases.

 \noindent Case 1: ($\delta=\gamma=0$).  In this case $\operatorname{Ad}^*({w,x,y,z})[\alpha,\beta,0,0]=[\alpha,\beta,0,0]$ for all $w,x,y,z$.  These are 1 point orbits determined by $\alpha$ and $\beta$.  The maximal subordinate algebra corresponding to any such orbit is the entire Lie algebra $\mathfrak{g}$, since $[A,B]\in \operatorname{Span}\{Y, Z\}$ for each $A,B\in \mathfrak{g}$.  Therefore $M_l=G$ and $G/M_l\cong 0$.  For any point $\{w,x,y,z\}\in G$, we write \[\{w, x, y, z\}=\exp\left(w, x, y-\frac{wx}{2}, z-\frac{x}{2}\left(y-\frac{wx}{2}+\frac{w^2}{6}\right)\right)\] and $\pi_{[\alpha, \beta, 0, 0]}$ is the one dimensional representation of $G$ given in $\mathbb{C}$ as
\begin{align*}
& \pi_{[\alpha,\beta,0,0]}\{w, x, y, z\}z
\\
& =e^{2\pi i [\alpha,\beta,0,0]\left(w, x, y-\frac{wx}{2},\,z-\frac{x}{2}\left(y-\frac{wx}{2}+\frac{w^2}{6}\right)\right)}z
\\
& =e^{2\pi i (\alpha w + \beta x)}z,
\end{align*}
for each $z\in \mathbb{C}$.

 \noindent Case 2:  ($\delta=0$, $\gamma\not=0$).  In this case $\operatorname{Ad}^*({w, x, y, z})[\alpha,\beta,\gamma,0]=[\alpha+x\gamma,\beta-w\gamma,\gamma,0],$ and so
 \[
 \operatorname{Ad}^*(G)[\alpha,\beta,\gamma,\delta]=\left\{[p, q, \gamma, 0]: p, q\in \mathbb{R}\right\}.
 \]
 These are 2-dimensional orbits parametrized by $\gamma$.  For any such orbit, the unitary irreducible representations induced by elements of the orbit are all unitarily equivalent and so it suffices to choose a convenient representative.  There is a one-to-one correspondence between the set
 \[
 R_2=\{[0, 0, \gamma,0]: \gamma \in \mathbb{R}^\times\}
  \]
 and the collection of orbits of this type.  Since $\gamma\not=0$,
 \[
 l_\gamma([W,X])=\gamma\not=0
 \]
 and so $\mathfrak{g}$ is not subordinate to $[0, 0, \gamma,0]$.  The three dimensional subalgebra $\mathfrak m=\operatorname{Span}\{X, Y, Z\}$ is Abelian and is therefore maximal subordinate to any element of $\mathfrak{g}^*$.  The subgroup

 \[
 M=\exp(\mathfrak m)=\{\{w, x, y, 0\}:w, x, y\in \mathbb{R}\}
 \]
 and $G/M\cong \mathbb{R}$. As indicated in \cite{TaylorBook}, $\pi_{[0, 0, \gamma, 0]}$ acts on

 \begin{align*}
\mathscr H_{\gamma}=\left\{ f: G \rightarrow  \mathbb{C}\Big{|}\right.&  f\in L^2(G/M)\text{ and }
\\
f(\exp(&q)g)=e^{2\pi i l_\gamma(q)}f(g)\text{ for each } q\in \mathfrak m \text{ and } g\in G\left.\phantom{\Big{|}}\right\}.
\end{align*}
First note that Haar measure $\mu$ on $G$ is given by $\mu(\exp(E))=\Lambda(E)$ where $\Lambda$ is Lebesgue measure on $\mathfrak{g}$, and so $\mathscr H_\pi:=L^2(G/M,\mu)\cong L^2(\mathbb{R},\Lambda)$.  We have that
\begin{align*}
\left(\pi_{[0, 0, \gamma, 0]}(\{w, x, y, z\})f\right)&(k)
\\
&= f(\{k,0,0,0\}*\{w,x,y,z\})
\\
&=f\left(\left\{k+w, x, y+kx, z+k \left(y+\frac{kx}{2}\right)\right\}\right)
\\
&=f\left(\left\{0, x, y+\frac{kx}{2}, z+k\left(y+\frac{kx}{2}\right)\right\}*\{k+w, 0, 0, 0\}\right)
\\
&=e^{2\pi i \gamma \left(y+\frac{kx}{2}\right)}f(k+w).
\end{align*}

\noindent Case 3:  ($\delta\not=0$).  We have that
\begin{align*}
\operatorname{Ad}^*(\{w, x, y, z\})[\alpha, \beta, \gamma, \delta] &
\\
&=\left[\alpha+x\gamma+(y-wx)\delta, \beta-w\gamma+\frac{w^2\delta}{2}, \gamma-w\delta, \delta\right].
\end{align*}
Defining $q=\gamma-w\delta$ we have that $w=\frac{\gamma-q}{\delta}$ and so
\begin{align*}
\operatorname{Ad}^*(\{w, x, y, z\})[\alpha, \beta, \gamma, \delta]&
\\
&=\left[\alpha+x\gamma+(y-wx)\delta, \left(\beta-\frac{\gamma^2}{2\delta}\right)+\frac{q^2}{2\delta}, q, \delta\right].
\end{align*}
Hence \[\operatorname{Ad}^*(G)[\alpha, \beta, \gamma, \delta]=\left\{\left[p, \left(\beta-\frac{\gamma^2}{2\delta}\right)+\frac{q^2}{2\delta}, q, \delta\right]: p, q\in \mathbb{R}\right\}.\]  These orbits are 2-dimensional parabolic cylinders parametrized by $\delta$ and the quantity $\beta-\frac{\gamma^2}{2\delta}$.  As in the previous case we have that
\[
R_3=\{[0,\beta,0,\delta]:\delta\in \mathbb{R}^\times,\beta\in\mathbb{R}\}
\]
is a collection of orbit representatives and $M=\operatorname{Span}\{X, Y, Z\}$ is a maximal subordinate subalgebra for each representative.  Therefore, $\mathscr H_{\beta,\delta}=L^2(\mathbb{R})$ and
\begin{align*}
\pi_{[0, \beta, 0, \delta]}(\{w, x, y, z\})f(k)&
\\
&=f\left(\{k, 0, 0, 0\}*\{w, x, y, z\}\right)\\
&=f\left(\left\{0, x, y+\frac{kx}{2}, z+k\left(y+\frac{kx}{2}\right)\right\}*\{k+w,0,0,0\}\right)\\
&=e^{2\pi i  \left(\beta x + \delta \left(z+k\left(y+\frac{kx}{2}\right)\right)\right)}f(k+w).
\end{align*}
\end{proof}

\section{The Weyl Functional Calculus for $G$}\label{s4}
 In Euclidean space, there is a well-developed theory of pseudo-differential operators and the corresponding symbolic calculus \cite{TaylorMBook1981}.  The classical Weyl functional calculus provides an analogous construction for the simplest step 2 nilpotent case.  A functional calculus for general nilpotent groups has been developed in \cite{BeltitaBeltita2009b}.  We will describe this functional calculus for $G$, and begin by stating the general construction for arbitrary nilpotent groups.
\begin{df} As above, let $K$ be an $n$ dimensional nilpotent Lie group with corresponding Lie algebra $\mathfrak k$.
\begin{enumerate}
\item Let $\xi_0\in \mathfrak k^*$ with corresponding co-adjoint orbit $\mathscr O$.  The {\it isotropy group of $K$ at $\xi_0$} is $K_{\xi_0}:=\{k\in K| \operatorname{Ad}^*(k)\xi_0=\xi_0\}$.
\item $K_{\xi_0}$ is a Lie group with corresponding {\it isotropy Lie algebra}

\[
\mathfrak k_{\xi_0}=\{X\in \mathfrak k|\xi_0\circ\operatorname{ad}(\mathfrak{k})X=0\}.
\]

\item Fix a sequence of ideals in $\mathfrak k$,\[\{0\}=\mathfrak k_0\subset\mathfrak k_1\subset\cdots\subset\mathfrak k_n=\mathfrak k\] such that $\dim(\mathfrak k_j/\mathfrak k_{j-1})=1$ and $[\mathfrak k,\mathfrak k_j]\subset \mathfrak k_{j-1}$ for $j=1, \dots, n$.  Pick any $X_j\in \mathfrak k_j\setminus \mathfrak k_{j-1}$ for $j=1, \dots, n$ so that the set $\{X_1, \dots, X_n\}$ is a {\it Jordan-H\"older basis} in $\mathfrak k$.
\item Consider the set of jump indices of the coadjoint orbit $\mathscr O$ with respect to the Jordan-H\"older basis,
\begin{eqnarray*}
J_{\xi_0}&=&\{j\in \{1,\dots,n\}|\mathfrak k_j\not\subseteq\mathfrak k_{j-1}+\mathfrak k_{\xi_0}\}\\
&=&\{j\in \{1, \dots, n\}|X_j\not\subseteq\mathfrak k_{j-1}+\mathfrak k_{\xi_0}\}
\end{eqnarray*}
 and then define the corresponding {\it predual of the coadjoint orbit $\mathscr O$},

 \[
 \mathfrak k_e:=\operatorname{Span}\{ X_j: j\in J_{\xi_0} \}.
 \]
 \item The {\it Fourier transform} $\mathscr S(\mathscr O)\rightarrow \mathscr S(\mathfrak{g}_e)$ is given by the formula
\[
\hat a(P)=\int_{\mathscr O}e^{-i\langle \xi,P\rangle}a(\xi)\operatorname{d} \xi \text{ for } P\in \mathfrak{g}_e,
\]
where $\operatorname{d} \xi $ is Liouville measure on $\mathscr O$.
 \item The {\it Weyl calculus} $\operatorname{Op}^\pi(\cdot)$ for the unitary representation $\pi$ is defined for every $a\in \mathscr S(\mathscr O)$ by
\[\operatorname{Op}^\pi(a)=\int_{\mathfrak k_e} \hat a(V) \pi (\exp_KV)\operatorname{d} V,\] where $\hat a(V)$ is the {\it Fourier transform} of $a\in \mathscr S(\mathscr O)$.    The operator $\operatorname{Op}^\pi(a)$ is called the {\it pseudo-differential operator} with symbol $a$.
\end{enumerate}
The following result appears in \cite{BeltitaBeltita2009b}.
\end{df}
\begin{thm}
The Weyl calculus $\operatorname{Op}^\pi$ has the following properties:
\begin{enumerate}
\item For every symbol $a\in \mathscr S(\mathscr O)$ we have $\operatorname{Op}^\pi(a)\in \mathscr B(\mathscr H)_\infty$ (the space of {\it smooth operators} for the representation $\pi$) and the mapping
\begin{align*}
\mathscr S(\mathscr O)\rightarrow \mathscr B(\mathscr H)_\infty &&a\mapsto\operatorname{Op}^\pi(a)
\end{align*}
 is a linear topological isomorphism.
\item For every $T\in \mathscr B(\mathscr H)_\infty$ we have $T=\operatorname{Op}^\pi(a)$ where $a\in \mathscr S(\mathscr O)$ satisfies the condition $\hat a(V)=\operatorname{Tr}(\pi(\exp_KV)^{-1}A)$ for every $V\in \mathfrak k_e$.
\end{enumerate}
\end{thm}
If $\pi$ is a representation of the nilpotent group $G$, then $\pi$ can be classified as in Proposition \ref{unirrepclassification}. If $\pi$ is of class 1 or class 2, then $\operatorname{Op}^\pi(\cdot)$ is understood \cite{TaylorBook}. From above results one can explicitly describe the Weyl functional calculus for class 3 representations of $G$.
\begin{prop}\label{Gweyl}
If $\pi$ is an irreducible unitary representation of $G$ of class 3 corresponding to the orbit $\mathscr O$ and $a\in \mathscr S(\mathscr O)$, then the Fourier transform of $a$ is given by

\[
\hat a\left( yY+wW \right)=\int_{\mathbb{R}^2}e^{-i(qy+pw)}a\left(q, p\right)\operatorname{d} q \operatorname{d} p
\]
and the pseudo-differential operator $\mathrm{Op}^\pi(a)$ is given for each $f\in L^2(\mathbb{R}^2)$ by
\[
\mathrm{Op}^\pi(a)f(k)=\int_{\mathbb{R}^2}\left[\int_{\mathbb{R}^2}e^{-i(qy+pw)}a(q,p)\operatorname{d} q\operatorname{d} p\right]e^{2\pi i(\delta ky+\frac{1}{2}\delta y w)}f(k+w)\operatorname{d} y \operatorname{d} w.
\]
\end{prop}
\begin{proof}
The basis $\{W,X,Y,Z\}$ is a Jordan-H\"older basis for $G$, and the predual of the co-adjoint orbit $\mathscr O$ is given by $\mathfrak{g}_e=\{W,Y\}$.  The chart $\mathscr O\rightarrow \mathbb{R}$
\begin{align*}
&pW^*+\left[\beta-\frac{q^2}{2\delta}\right]X^*+qY^*+\delta Z^*\mapsto\left(p,q\right)
\end{align*}
is a map which brings Liouville measure on $\mathscr O$ to Lebesgue measure on $\mathbb{R}$.  Direct substitution implies that the Fourier transform is given by \[\hat a\left(yY+wW\right)=\int_{\mathbb{R}^2}e^{-i(qy+pw)}a\left(q,p\right)\operatorname{d} q \operatorname{d} p.\]  For $\left(\pi(\{w, x, y, z\})f\right)(k)=e^{2\pi i(\beta x + \delta(z+k(y+\frac{kx}{2})))}f(k+w)$ and $(w, 0, y, 0)\in \mathfrak g_e$ we have that \[\pi(\exp(w, 0, y, 0))f(k)=\pi(\{w, 0, y, 0\})f(k)=e^{2\pi i(\delta(ky+\frac{k^2x}{2}))}f(k+w),\] and direct substitution yields the result.

\end{proof}

\section{L\'{e}vy Processes in $G$} \label{s5}
 The expository material of this section can be found in \cite{LiaoBook}.  Suppose that $K$ is an arbitrary (not necessarily nilpotent) Lie group with Lie algebra $\mathfrak k$.  A L\'evy process in $K$ is a $K$-valued stochastic process $X_{t},t\geqslant 0$ which satisfies the following
\begin{enumerate}
\item $X_{t}$ has stationary and independent left increments, where the increment between $s$ and $t$ with $s\leqslant  t$ is $X_{t}(s)^{-1}X_{t}$.
\item $X_{t}(0)=e$ a.s.
\item $X_{t}$ is stochastically continuous, i.e.
\[
\lim_{s\rightarrow t}P(X_{s}^{-1}X_{t}\in A)=0
 \]
for all $A\in \mathcal B(K)$ such that $e\not\in\overline A$.
\end{enumerate}
 Let $C_0(K)$ be the Banach space (with respect to the supremum norm) of functions on $K$ which vanish at infinity.  Just as in the Euclidean case, one obtains a Feller semigroup on $C_0(K)$ by the prescription
 \[
(T(t)f)(k)=\mathbf E(f(k X_{t})),
 \]
 for each $t\geqslant  0$, $k \in K$, $f\in C_0(K)$ and its infinitesimal generator will be denoted as $\mathcal L$.

We fix a basis $\{Z_1, \dots, Z_n\}$ for $\mathfrak k$ and define a dense subspace $C_2(K)$ of $C_0(K)$ as follows:
\begin{align*}
& C_2(K)=\left\{f\in C_0(K):\right.
\\
& \left. Z^L_i(f)\in C_0(K)\textrm{ and }Z^L_iZ^L_j(f)\in C_0(K)\textrm{ for all }1\leqslant  i,j\leqslant  n \right\},
\end{align*}
where $Z^L$ denotes the left invariant vector field associated to $Z\in \mathfrak k$.

 In \cite{Hunt1956}, Hunt proved that there exist local coordinate functions $y_i\in C_2(K)$, $1\leqslant  i\leqslant  n$ so that each \[y_i(e)=0\textrm{   and   }Z^L_iy_j(e)=\delta_{ij},\] and a map $h\in \operatorname{Dom}(\mathcal L)$ which is such that:
\begin{enumerate}
\item $h>0$ on $K-\{e\}$.
\item There exists a compact neighborhood of the identity $U$ such that for all $\tau\in U$, \[h(\tau)=\sum_{i=1}^ny_i(\tau)^2.\]  Any such function is called a {\it Hunt function} in $K$.  A positive measure $\nu$ defined on $\mathcal B(Q-\{e\})$ is called a {\it L\'evy measure} whenever
\begin{equation}\label{HuntFunction}
\int_{Q-\{e\}}h(\sigma)\nu(\operatorname{d} \sigma)<\infty.
\end{equation}
\end{enumerate}
\begin{thm}[Hunt]
Let $X_{t}$ be a L\'evy process in $K$ with infinitesimal generator $\mathcal L$ then,
\begin{enumerate}
\item $C_2(K)\subset \operatorname{Dom}(\mathcal L)$.
\item For each $\tau\in K$, $f\in C_2(K)$
\begin{align}
\mathcal L(\tau)=&\sum_{i=1}^nb_iZ_i^Lf(\tau)+\sum_{i,j=1}^nc_{ij}Z^L_iZ_j^Lf(\tau)\nonumber\\
&+\int_{K-\{e\}}(f(\tau\sigma)-f(\tau)-\sum_{i=1}^ny_i(\sigma)Z_i^Lf(\tau))\nu(\operatorname{d} \sigma),\label{applebaumhunt}
\end{align}
where $b=(b_1,\dots,b_n)\in \mathbb{R}^n$, $c=(c_{ij})$ is a non-negative-definite, symmetric $n\times n$ real-valued matrix and $\nu$ is a L\'evy measure on $K-\{e\}$.
\end{enumerate}
Furthermore, any linear operator with a representation as in \ref{applebaumhunt} is the restriction to $C_2(K)$ of a unique weakly continuous, convolution semigroup of probability measures in $K$.
\end{thm}
 Let $\mathscr H$ be a complex, separable Hilbert space and $U(\mathscr H)$ be the group of unitary operators in $\mathscr H$.  Let $\pi:K\rightarrow U(\mathscr H)$ be a strongly continuous unitary representation of $K$ in $\mathscr H$ and let $C^\infty(\pi)=\{\psi\in \mathscr H;k\rightarrow \pi(k)\psi\textrm{ is }C^\infty\}$ be the dense linear space of {\it smooth vectors for $\pi$} in $\mathscr H$.  Define a strongly continuous contraction semigroup $\mathcal T_t$ of linear operators on $\mathscr H$ by
 \[
 \mathcal T_t\psi=\mathbf E(\pi(X_{t})\psi)
 \]
 for each $\psi\in \mathscr H$. Let $L^{\pi}$ denote the infinitesimal generator of this semigroup. It follows from the work in \cite{ApplebaumCohen2004} that $C^\infty(\pi)\subseteq\operatorname{Dom}(\mathcal L^\pi)$ and for $f\in C^\infty(\pi)$ we have
 \begin{align}
 \mathcal L^{\pi} f=&\sum_{i=1}^n b_i\operatorname{d} \pi(Z_i)f+\sum_{i,j=1}^nc_{ij}\operatorname{d} \pi(Z_i)\operatorname{d} \pi(Z_j)f+\nonumber\\
&+\int_{K-\{e\}}\left(\pi(\sigma)-I-\sum_{i=1}^ny_i(\sigma)\operatorname{d} \pi(Z_i)\right)f\nu(\operatorname{d} \sigma).\label{acgenerator}
\end{align}
We now investigate $\mathcal L^\pi$ where $K=G$.  Since $G$ is nilpotent, the Haar measure $\operatorname{d} \sigma$ is related to Lebesgue measure on $\mathfrak g$ via the exponential map.  Therefore it will be convenient to adopt exponential coordinates in $G$.  To this end we impose the identification of $(w, x, y, z)$ with $\exp( w, ,y , z )$.  Fix real numbers $\beta$ and $\delta\not=0$.   Let $\pi=\pi_{\delta, \beta}$ be a representation of class 3.  Define
 \begin{align*}
 & Kf(k)=kf(k),
 \\
 & Df(k)=\frac{1}{i}\frac{\operatorname{d} f}{\operatorname{d} k}.
 \end{align*}
We have that
 \begin{align}
 \left(\pi( w, x, y, z )f\right)(k)=e^{2\pi i\left(\left(\beta x+\delta \left(z+\frac{xy}{2}+\frac{w^2x}{6}\right)\right)I+\left(y+\frac{wx}{2}\right)K+\frac{x}{2}K^2\right)+2\pi i wD}f(k)\label{thirdtyperep}
 \end{align}
 and
\begin{align*}
&\operatorname{d}\pi (W)=2\pi i D,
\\
&\operatorname{d}\pi(X)=2\pi i \beta I+\pi i K^2,
\\
&\operatorname{d}\pi (Y)=2\pi i \delta K,
\\
&\operatorname{d}\pi (Z)=2\pi i \delta I.
\end{align*}
Denote
\begin{align*}
& \mathcal L^\pi_1:=\sum_{i=1}^n b_i\operatorname{d} \pi(Z_i),
\\
& \mathcal L^\pi_2:=\sum_{i,j=1}^nc_{ij}\operatorname{d} \pi(Z_i)\operatorname{d} \pi(Z_j),
\\
& \mathcal L^\pi_3:=\int_{G-\{e\}}\left(\pi(\sigma)-I-\sum_{i=1}^ny_i(\sigma)\operatorname{d} \pi(Z_i)\right)f\nu(\operatorname{d} \sigma).
\end{align*}
Then the drift part can be written as follows.
\begin{align}
\mathcal L^\pi _1&= b_1(2\pi i\delta I)+b_2(2\pi i\delta K)+b_3(2\pi i\beta I+\pi i K^2)+b_4(2\pi iD).\label{Lpidrift}
\end{align}
Using the Weyl functional calculus described in Proposition \ref{Gweyl}, $\mathcal L^\pi_1$ is a pseudo-differential operator with symbol given by
\begin{align*}
\mathcal S^\pi _1&= 2\pi i\delta b_1+2\pi i\delta b_2 t+b_3(2\pi i\beta +\pi i t^2)+2\pi ib_4\frac{\partial}{\partial t}.
\end{align*}
The Brownian part can be expressed
\begin{align}
\mathcal L^\pi _2=   &c_{11}(-4\pi ^2 \delta ^2 I)+c_{22}(-4\pi^2\delta ^2 K^2)\label{Lpibrownian}\\
&+c_{33}(-4\pi^2\beta^2 I - 4\pi^2 \beta K^2-\pi^2 K^4)\notag\\
&+c_{44}(-4\pi^2 D^2)+ 2c_{12}(-2\pi ^2\delta^2K)+2c_{13}(-4\pi^2\delta\beta I-2\pi^2\delta K^2)\notag\\
&+2c_{14}(-4\pi^2\delta D)+2c_{23}(-4\pi^2\delta\beta K-2\pi^2\delta K^3)+c_{24}(-4\pi^2\delta KD)\notag\\
&+c_{34}(-4\pi^2\beta D-2\pi^2K^2D)+c_{42}(-4\pi^2\delta(KD+I))\notag\\
&+c_{43}(-4\pi^2\beta D-2\pi^2(2K+K^2D)),\nonumber
\end{align}
which is a pseudo-differential operator with symbol
\begin{align*}
\mathcal S^\pi _2=   &-4\pi ^2 \delta ^2 c_{11}-4\pi^2\delta ^2c_{22} t^2+c_{33}(-4\pi^2\beta^2  - 4\pi^2 \beta t^2-\pi^2 t^4)\\
&+ c_{44}\left(-4\pi^2 \frac{\partial^2}{\partial t^2}\right)-4\pi ^2\delta^2c_{12}t+2c_{13}(-4\pi^2\delta\beta -2\pi^2\delta t^2)\\
&- 8\pi^2\delta c_{14} \frac{\partial}{\partial t}+2c_{23}(-4\pi^2\delta\beta t-2\pi^2\delta t^3)-8\pi^2\delta c_{24}t\frac{\partial}{\partial t}\\
&+ c_{34}\left(-4\pi^2\beta \frac{\partial}{\partial t}-2\pi^2t^2\frac{\partial}{\partial t}\right)-4\pi^2\delta c_{42}\left(t\frac{\partial}{\partial t}+1\right)\\
&+ c_{43}\left(-4\pi^2\beta \frac{\partial}{\partial t}-2\pi^2\left(2t+t^2\frac{\partial}{\partial t}\right)\right).
\end{align*}
 Before expressing the jump part $\mathcal L^\pi_3$, observe that (\ref{thirdtyperep}) can be rewritten as
\begin{align*}
\pi( w, x, y, z )f(k)=\exp(i\Phi( w, x, y, z ))f(k)\\
\end{align*}
where
\begin{align*}
\Phi( w, x, y, z )=  & 2\pi\left[  \delta I\left(z+\frac{1}{2}w^2x+\frac{1}{2}wy\right)+  \delta K\left(y+\frac{1}{2}wx\right)\right.\\
&\left.\phantom{\frac{1}{2}}+\left( \beta I+\frac{1}{2}  K^2\right)(x)+  wD\right]
\end{align*}
is essentially self-adjoint.  This form suggests the following choices for local coordinate functions
\begin{align*}
& y_1( w, x, y, z)=w\chi_B( w, x, y, z ),
\\
& y_2( w, x, y, z )=x\chi_B( w, x, y, z ),
\\
& y_3( w, x, y, z )=\left(y+\frac{1}{2}wx\right)\chi_B( w, x, y, z ),
\\
& y_4( w, x, y, z)=\left(z+\frac{1}{2}w^2x+\frac{1}{2}wy\right)\chi_B( w, x, y, z ),
\end{align*}
where $y_i( w, x, y, z )=y_i(\exp( w, x, y, z ))$, $B=\exp(B( 0, 1 ))$.  With respect to these local coordinate functions we have that
\begin{align*}
\mathcal L^\pi_3  =\\
\int_{\mathbb{R}^4\backslash\{0\}}&\left(\pi( w, x, y, z)-I- i\Phi( w, x, y, z) \chi_B( w, x, y, z)  \right)\nu(\operatorname{d}z \operatorname{d}y \operatorname{d}x \operatorname{d}w).
\end{align*}
Therefore $\mathcal L^\pi_3$ is a pseudo-differential operator with symbol
\begin{align*}
\mathcal S^\pi_3  =
\\
\int_{\mathbb{R}^4\backslash\{ 0 \}}&\left(\tau( w, x, y, z )-I-i\Theta( w, x, y, z)\chi_B( w, x, y, z)    \right)\nu(\operatorname{d}z \operatorname{d}y \operatorname{d}x \operatorname{d}w),
\end{align*}
where
\begin{align*}
\tau( w, x, y, z )= &\exp(i\Theta( w, x, y, z ))
\end{align*}
for
\begin{align*}
\Theta( w, x, y, z )= &2\pi\left[  \delta \left(z+\frac{1}{2}w^2x+\frac{1}{2}wy\right)+  \delta t\left(y+\frac{1}{2}wx\right)\right.\\
&\left.+\left( \beta +\frac{1}{2}  t^2\right)(x)+  w\frac{\partial }{\partial t}\right]
\end{align*}
and $\pi$ is as in (\ref{thirdtyperep}).
We are now ready to state the main theorem of this paper.
\begin{thm}
The operator $\mathcal L^\pi$ is a pseudo-differential operator.  Moreover, the space $C^\infty_c(\mathbb{R})$ is a core for $\mathcal L^\pi$.
\end{thm}
\begin{proof}
 We have that \[\mathcal L^\pi=\mathcal L^\pi_1+\mathcal L^\pi_2+\mathcal L^\pi_3,\] and consequently we have shown that $\mathcal L^\pi$ is pseudo-differential with symbol \[\mathcal S^\pi=\mathcal S^\pi_1+\mathcal S^\pi_2+\mathcal S^\pi_3.\]  We write $\mathcal L^\pi_3=\mathcal L^\pi_{3,1}+\mathcal L^\pi_{3,2}$ with
\begin{align*}
\mathcal L^\pi_{3,1}= &\int_{B^c}\left(\pi( w, x, y, z )-I\right)\nu(\operatorname{d} z \operatorname{d} y \operatorname{d} x \operatorname{d} w)\\
\mathcal L^\pi_{3,2}= &\int_{B-\{0\}}\left(\pi( w, x, y, z )-I\phantom{\frac{1}{2}}-i\Phi( w, x, y, z )\right)\nu(\operatorname{d} z  \operatorname{d} y  \operatorname{d} x  \operatorname{d} w).\\
\end{align*}
For each $f\in C^\infty_c(\mathbb{R})$, we have that
\begin{align*}\|\mathcal L^\pi _{3,1}f\|\leqslant  &\int_{B^c}\left\|\left(\pi( w, x, y, z )-I\right)f\right\|\nu(\operatorname{d} z \operatorname{d} y \operatorname{d} x \operatorname{d} w)\\
\leqslant &2\nu(B^c)\|f\|.
\end{align*}
Let $P( w, x, y, z )$ denote the projection-valued measure associated to the spectral decomposition of the self adjoint operator $\Phi$.  By the spectral theorem and Taylor's theorem, and referring again to (\ref{thirdtyperep}) we see that
\begin{align*}
\|(\pi( w, x, y, z )-I-&i\Phi( w, x, y, z ))f   \|^2\\
= &\int_{\mathbb{R}^4} \left|e^{i\lambda}-1-i\lambda\right|^2\|P( w, x, y, z )(\operatorname{d} \lambda)f\|^2\\
\leqslant  &\frac{1}{4}\int_{\mathbb{R}^4} |\lambda|^4\|P( w, x, y, z )(\operatorname{d} \lambda)f\|^2\\
= &\frac{1}{4}\left\|2\pi\left[ \delta I\left(z+\frac{1}{2}w^2x+\frac{1}{2}wy\right)+ \delta K\left(y+\frac{1}{2}wx\right)\right.\right.\\
&\left.\left.+\left(\beta I+\frac{1}{2} K^2\right)(x)+ wD\right]^2f\right\|^2\\
\leqslant  & \pi^2\left\|\left[y_1( w, x, y, z )D+y_2( w, x, y, z )\left(\beta I+\frac{1}{2}K^2\right)\right.\right.\\
&\left.\left.\phantom{\frac{1}{2}}+y_3( w, x, y, z )\delta K+y_4( w, x, y, z )\delta I\right]^2f\right\|^2\\
\leqslant  &16\pi^2 C_f^2\cdot h^2( w, x, y, z ).
\end{align*}
The last inequality follows from Young's inequality. The Hunt function $h$ corresponds to the local coordinate functions $\{y_i\}_{i=1}^4$ and
\begin{align*}
C_f=&((\beta+\delta)^2+\delta)\|f\|+(2\delta(\beta+\delta)+1)\|Kf\|+(\delta^2+\beta+\delta)\|K^2f\|\\
&+\delta\|K^3f\|+\frac{1}{4}\|K^4f\|+2(\beta+\delta)\|Df\|+2\delta\|KDf\|+2\|K^2Df\|\\
&+\|D^2f\|.
\end{align*}
Therefore we have that
\begin{align*}
\|\mathcal L^\pi_{3,2}f\|\leqslant & 4\pi C_f \int_{B}h( w, x, y, z )\nu(\operatorname{d}w \operatorname{d}x \operatorname{d}y \operatorname{d}z),
\end{align*}
and the latter integral is finite by \eqref{HuntFunction}.  Applying these bounds for $\mathcal L_{3,1}^\pi$ and $\mathcal L_{3,2}^\pi$ and the expressions (\ref{Lpidrift}) and (\ref{Lpibrownian}) there exist non-negative constants $\omega(f)_{ij}$ such that
\begin{align}
\|\mathcal L^\pi f\|\leqslant  \sum_{i=1}^4\sum_{j=1}^2\omega(f)_{ij}\|K^iD^jf\|.\label{Lpibound}
\end{align}
Let $f\in \textrm{Dom}(\mathcal L^\pi)$, then we can find $(f_n,n\in \mathbb{N})$ in $C^\infty_c(\mathbb{R})$ such that
\[
\lim_{n\rightarrow \infty}\|f_n-f\|= 0.
\]
Applying \eqref{Lpibound} to the sequence $f_n-f_m$, we deduce by integration by parts and the Schwarz inequality that $\lim_{m,n\rightarrow\infty}\Vert\mathcal L^\pi(f_n-f_m)\Vert=0$.  Hence the sequence $(\mathcal L^\pi f_n, n\in \mathbb{N})$ is Cauchy and so convergent to some $g\in L^2(\mathbb{R})$.  The operator $\mathcal L^\pi$ is closed, hence $g=\mathcal L^\pi f$ and the result is established.
\end{proof}

\begin{example}
Let $B^1_t$ and $B^2_t$ be independent one-dimensional Brownian motions and define $X_t=\exp(B^1_t,B^2_t,0,0)$. Then its generator is the sub-Laplacian
\begin{align*}
\mathcal L=\left(\frac{\partial}{\partial w}\right)^2+\left(\frac{\partial}{\partial x}+w\frac{\partial}{\partial y}+\frac{w^2}{2}\frac{\partial}{\partial z}\right)^2.
\end{align*}
If $\pi=\pi_{\delta,\beta}$ is a representation of class 3, the quantization of this generator is given as
\begin{align*}
\mathcal L^\pi=-\pi^2\left(4D^2+4\beta^2I+4\beta K^2 +K^4\right).
\end{align*}
The symbol of $\mathcal L^\pi$ is the operator $S^\pi$ defined for each $\varphi(t)\in C^\infty_c(\mathbb{R})$ by
\begin{align*}
\mathcal S^\pi\varphi(t)=-\pi^2\left(4\frac{\partial^2\varphi}{\partial t^2}+4\beta^2+4\beta t^2\varphi(t) +t^4\varphi(t)\right).
\end{align*}
In \cite{baudoincoutin} it is shown that fractional Brownian motion in Carnot groups exhibits a scaling property reminiscent of the property for Brownian motion in $\mathbb{R}^n$.  If one defines
\begin{align*}
V_1&=\operatorname{Span}\{ W, X \},\\
V_2&=\operatorname{Span}\{Y\},\\
V_3&=\operatorname{Span}\{Z\},\\
\end{align*}
then $\mathfrak{g}=V_1\oplus V_2\oplus V_3$ and it is clear that $G$ has the structure of a Carnot group.  A fractional Brownian motion is not a L\'evy process unless o the Hurst parameter $H$ is equal to $\frac{1}{2}$.  By a standard application of Ito's Lemma
\begin{align*}
X_t= & \left\{B^1_t, B^2_t, \frac{1}{2}\int_0^t(B^1_sdB^2_s-B^2_sdB^1_s)+\frac{1}{2}B^1_tB^2_t, \right.
\\
&\left. \int_0^t\left(\int_0^s \left(\frac{1}{6}B^2_rdB^1_r-\frac{1}{3}B^1_rdB^2_r\right)\right)dB^1_s+\right.\\
&\left.\frac{B^1_t}{4}\int_0^t(B^1_sdB^2_s-B^2_sdB^1_s)+\frac{(B^1_t)^2B^2_t}{6}\right\}.
\end{align*}
It is easy to verify that $X_t$ solves the stochastic differential equation
\begin{align*}
& dX_t=X_t(WdB^1_t+XdB^2_t),
\\
& X_0=e.
\end{align*}
Finally, the scaling property as formulated in \cite[Proposition 3.8]{baudoincoutin} implies that
\begin{align*}
(X_{ct})_{t\geqslant  0}{\buildrel law \over =} (\Delta_{\sqrt c}X_t)_{t\geqslant  0},
\end{align*}
where
\begin{align*}
\Delta_{\sqrt c}\{w, x, y, z\}=\{\sqrt c w, \sqrt cx, cy, \sqrt{c^3}z\}.
\end{align*}
\end{example}

\bibliographystyle{amsplain}
\bibliography{GordinaHaga}

\end{document}

%% file: GordinaHaga-title.tex
\def\enddoc{\end{document}}

\title[L\'evy Processes in a Step 3 Nilpotent Lie Group]{L\'evy Processes in a Step 3 Nilpotent Lie Group}

\author[Gordina]{Maria Gordina{$^{\dagger}$}}
\thanks{\footnotemark {$\dagger$} This research was supported in part by NSF Grant DMS-1007496.}
\address{Department of Mathematics\\
University of Connecticut\\
Storrs, CT 06269, USA} \email{maria.gordina@uconn.edu}
\author[Haga]{John Haga{$^{\ast}$}}
\thanks{\footnotemark {$^{\ast}$} This research was supported  in part by NSF Grant DMS-1007496.}
\address{Department of Mathematics\\
University of Connecticut\\
Storrs, CT 06269, USA} \email{john.haga@uconn.edu}

\keywords{L\'{e}vy processes, nilpotent group, heat kernel, pseudo-differential operator}
\subjclass{Primary 60G51, 35S05, 58J65; Secondary 60J25}


\date{\today \ \emph{File:\jobname{.tex}}}

\begin{abstract}

The infinitesimal generators of L\'evy processes in Euclidean space are pseudo-differential operators with symbols given by the L\'evy-Khintchine formula. This classical analysis relies heavily on Fourier analysis which in the case when the state space is a Lie group becomes much more subtle.  Still the notion of pseudo-differential operators can be extended to connected, simply connected nilpotent Lie groups by employing the Weyl functional calculus.  With respect to this definition, the generators of L\'evy processes in the simplest step 3 nilpotent Lie group $G$ are pseudo-differential operators which admit $C_c(G)$ as its core.

\end{abstract}

\maketitle

\tableofcontents

\renewcommand{\contentsname}{Table of Contents}